# Global existence of solutions to a tear film model with locally elevated evaporation rates

Yuan Gao[a], Hangjie Ji[b], Jian-Guo Liu[b,c], Thomas P. Witelski[b]

[a]*School of Mathematical Sciences, Fudan University, P.R. China*
[b]*Department of Mathematics, Duke University*
[c]*Department of Physics, Duke University*

January 3, 2017

**Abstract**

Motivated by a model proposed by Peng et al. [*Advances in Coll. and Interf. Sci.* 206 (2014)] for break-up of tear films on human eyes, we study the dynamics of a generalized thin film model. The governing equations form a fourth-order coupled system of nonlinear parabolic PDE for the film thickness and salt concentration subject to non-conservative effects representing evaporation. We analytically prove the global existence of solutions to this model with mobility exponents in several different ranges and the results are then validated against PDE simulations. We also numerically capture other interesting dynamics of the model, including finite-time rupture-shock phenomenon due to the instabilities caused by locally elevated evaporation rates, convergence to equilibrium and infinite-time thinning.

*Keywords:* global existence, tear film, fourth-order nonlinear partial differential equations, rupture, thin film equation, evaporation, osmolarity, finite-time singularity

## 1. Introduction

In this article, we study the regularity of solutions to a one-dimensional nonlinear partial differential equation system for a fluid film height $h(x,t)$ and salt concentration (also called the osmolarity) $s(x,t)$ on a finite domain, $0 \leq x \leq L$,

$$h_t = -(h^n h_{xxx})_x - h^m(\bar{S} - s), \tag{1.1a}$$

$$s_t = s_{xx} + \left(\frac{h_x}{h} - h^{n-1} h_{xxx}\right) s_x + s(\bar{S} - s) h^{m-1}. \tag{1.1b}$$

This family of PDEs is motivated by a non-conservative lubrication model for evaporating tear films on human eyes. Based on the model proposed by Peng et al. [20], a spatial variation in a thin lipid layer on the tear film leads to locally elevated evaporation rates of the tear film, which in turn affects the local salt concentration in the liquid film. In our model (1.1) the influences of the lipid layer thickness on osmolarity are included in the effective salt capacity function, $\bar{S}(x) \in L^\infty([0, L])$. This will be taken to be a given positive function with increased values over some portion of the domain, corresponding to elevated evaporation rates (and decreased lipid concentrations). Details of the formulations from the physical model will be discussed further in section 2.

The mobility exponents $n$ and $m$ in (1.1) are introduced to analyze and separate the influences of the conservative and non-conservative fluxes in the model respectively. Starting from initial data $(h_0(x), s_0(x))$ at time $t = 0$ which satisfy $h_0 > 0$ and $0 \leq s_0 \leq \|\bar{S}\|_\infty$, the dynamics will be subject to no-flux and normal-contact boundary conditions

$$h_x(0) = h_x(L) = 0, \quad h_{xxx}(0) = h_{xxx}(L) = 0, \quad s_x(0) = s_x(L) = 0. \tag{1.2}$$

*Email address:* `hangjie@math.duke.edu` (Hangjie Ji)

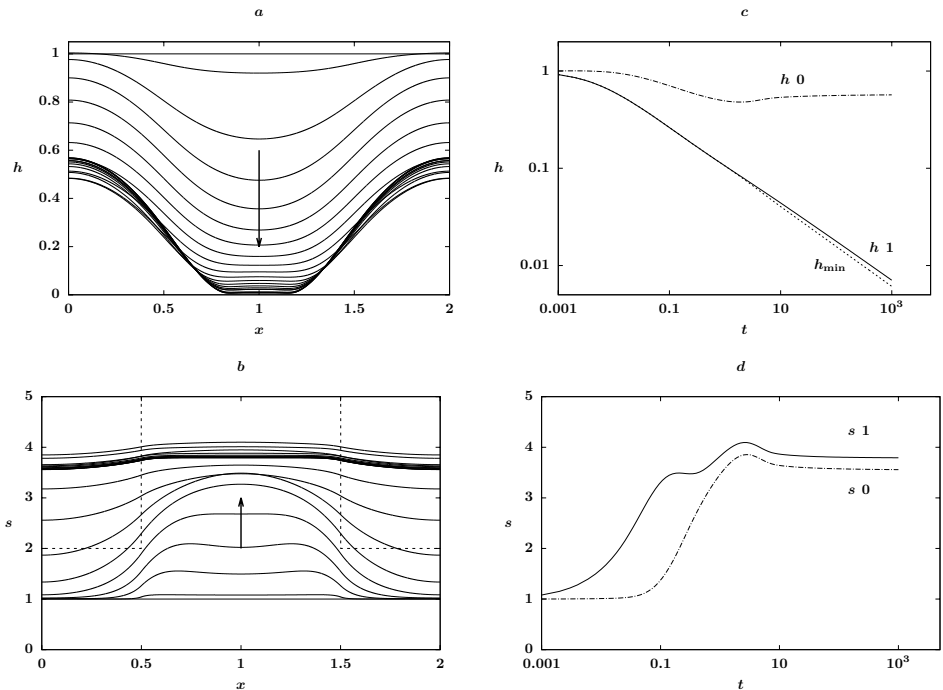

Figure 1: Numerical solution of (1.1) with $(m,n) = (3.5, 4.5)$ starting from constant initial data $h_0 = s_0 = 1$ for $0 \leq x \leq 2$ driven by non-conservative flux with $\bar{S}(x)$ given by (1.5): $(a)$ and $(b)$ are solution profiles for $h$ and $s$, where the capacity $\bar{S}$ is plotted in dashed line; $(c)$ and $(d)$ represent the evolution of solutions at $x = 0$, $x = 1$, respectively. The minimum film thickness $h_{\min}(t)$ monotonically decreases and approaches 0 as $t \to \infty$, while the osmolarity $s$ is locally elevated in the early stage and the spatial variations decay in the later stage.

The total mass of salt, $Q$, in the liquid film is of fundamental interest, and one can obtain the conservation of this quantity by multiplying (1.1a) by $s$, multiplying (1.1b) by $h$ and then integrating the sum of the two equations,

$$Q(t) = \int_0^L hs \, dx, \qquad Q(t) = Q(0) = \int_0^L h_0 s_0 \, dx. \tag{1.3}$$

We will use the conservation of $Q$ in the proof of the boundedness of the $h$ and the study of equilibrium solutions of the system. The total mass of the liquid film in the system is not conserved and its rate of change can be obtained from integrating (1.1a) and applying the boundary conditions.

$$M(t) = \int_0^L h \, dx, \qquad \frac{dM}{dt} = -\int_0^L h^m (\bar{S} - s) \, dx. \tag{1.4}$$

We note that with initial data $0 < s_0 < \|\bar{S}\|_\infty$ and $h_0 > 0$, the osmolarity is always bounded by $\|\bar{S}\|_\infty$. This result will be shown formally in the proof of Theorem 1 in section (3) by the weak form of maximum principle. Namely, $s(x,t)$ is not guaranteed to be bounded by $\bar{S}(x)$ pointwise; consequently it is not clear if $M(t)$ is necessarily decreasing in time.

The imposed osmolar capacity $\bar{S}(x)$ is essential in determining the dynamics of the model (1.1). Although in previous models [20], it was assumed that $\bar{S}$ is smooth and takes the form of a Gaussian distribution, the regularity of solutions to PDEs (1.1) are not sensitive to the smoothness of $\bar{S}$. For instance, with $\bar{S}$ given by a positive step function,

$$\bar{S}(x) = \begin{cases} 100 & \text{for } 0.5 < x < 1.5, \\ 2 & \text{otherwise,} \end{cases} \tag{1.5}$$

a typical evolution of $h$ and $s$ over the domain $0 \leq x \leq 2$ with $(m,n) = (3.5, 4.5)$ is presented in Figure 1, with figures $(a)$ and $(b)$ showing the dynamics of height and salt concentration profiles in (1.1) and figures $(c)$ and (d) illustrating the evolution of some key properties for those profiles. Spatially uniform initial conditions for both thickness and salt concentration are used, corresponding to states produced by opened eyes after a blink [4]. Starting from normalized initial conditions $h_0 = s_0 = 1$, the film thickness $h$ decreases with an increasing spatial variation driven by the locally elevated $\bar{S}$, while the osmolarity rises due to the evaporation and is elevated more quickly near the center of the domain. Later, the symmetry of the film is broken and $h$ reaches its minimum at $x \approx 0.8$ and $x \approx 1.2$. This is difficult to see in plot (a), but is depicted in the plot (c) where $h_{\min}(t) = \min_x h(x,t)$ is no longer attained at the center of the domain. The localized elevation in osmolarity $s$ is a transient effect with $h_{\min}(t)$ approaching zero algebraically and $s$ evolving slowly as $t \to \infty$.

The mobility exponents, $m$ and $n$, are crucial in determining the qualitative behavior of the PDE model. In section 2 we will review the motivating tear film model which corresponds to (1.1) with exponents $(m,n) = (0,3)$. The modified tear film equations (1.1) are then proposed. In section 3, the global and local existence of solutions to (1.1) in several different ranges of $(m,n)$ are proved analytically. These regularity results, together with other interesting dynamics of the model, are then further investigated with numerical simulations in section 4.

## 2. A physical model for tear film break-up driven by evaporation

Human eyes are coated with a thin layer of precorneal tear film, which is a bio-fluid with a complex composition, but which can be approximated in terms of a viscous fluid with dissolved salt and a lipid layer. Tear film thinning and break-up during interblink periods are observable clinically and play a key role in dry eye disorders. While many mechanisms have been proposed in the literature [11, 18, 23], it is now generally agreed that evaporation is one of the most important factors for the tear film break-up phenomenon [16]. In addition to evaporation effects, capillarity and osmolarity also contribute to the dynamics of the tear film. In particular, experiments have shown that the local increase of osmolarity is evident along with the reduced tear film thickness as the break-up occurs. For a thorough discussion on the dynamics of tear film, we refer to [4].

Since the average aqueous tear film thickness (approximately $10^{-6}$ m) [15] is much thinner than the average radius of curvature of the eye (about $10^{-2}$ m), we follow the literature [4, 3] to assume that the substrate underneath the tear film is flat. In addition, it is appropriate to use the lubrication approximation to model the evolution of the tear film, since the exposed surface of the eye (called the palpebral fissure) is about $10^3$ times larger than the average tear film thickness. The dynamics of thin viscous film flows have been studied extensively for the past decade due to their fundamental importance in coating flows, painting, biological applications like tear films [4], and other applications in science and engineering. There is a large literature on the numerical and modeling studies of evaporating thin films [9], and some of the results can be applied to the study of tear films. For the characterization of thin film flows with surfactant [26, 21], a suspension of heavy particles [8], or drying paint layers [12, 10], an additional PDE for surfactant concentration or particle concentration is usually incorporated into the PDE system.

The conservation of water leads to the dimensionless governing equation for the film thickness $h$,

$$\frac{\partial h}{\partial t} = -\frac{\partial}{\partial x}(uh) - J_e + J_w, \tag{2.1}$$

with the flow velocity $u$ given by [20]

$$u = -h^2 \frac{\partial p}{\partial x}, \tag{2.2}$$

where the dynamical pressure $p$ is given by the combination of the generalized conjoining pressure $\Pi(h)$ and the linearized curvature $h_{xx}$,

$$p = \Pi(h) - h_{xx}, \tag{2.3}$$

and $t$ and $x$ are temporal and spatial variables. Here we take the conjoining pressure $\Pi(h) = -\epsilon/h^3$ with the rescaled Hamaker constant $\epsilon > 0$ which represents the wetting property of the corneal surface. The two non-conservative contributions $J_e$ and $J_w$ correspond to the evaporative flux and osmotic weeping flux and will be described below in detail.

The dynamics of the salt concentration is governed by a second dimensionless evolution equation describes the local conservation of salt in terms of diffusion and convective transport of the aqueous film,

$$\frac{\partial}{\partial t}(hs) = \frac{\partial}{\partial x}\left(h\frac{\partial s}{\partial x} - ush\right). \tag{2.4}$$

This form guarantees that the total salt mass is conserved

$$\int hs\ dx = \int h_0 s_0\ dx. \tag{2.5}$$

The osmotic weeping flux, $J_w$, is assumed to be proportional to the difference between the salt concentration in the tear film and a reference osmolarity constant, $S_0$,

$$J_w = s - S_0. \tag{2.6}$$

Several different mathematical models have been studied [17, 24, 5, 25] for the mechanisms of the tear film break-up driven by evaporation, and the major differences among these models are in the physical interpretation of the evaporation effects involved in the tear film. For instance, Braun [4] derived the form of evaporative mass flux

$$J_e = \frac{E(1+\delta p)}{\bar{K}+h}, \tag{2.7}$$

where $\bar{K} > 0$ measures the thermal resistances to mass transfer at the fluid-vapor interface, and $E > 0$ characterizes the ratio of viscous timescale to the non-conservative timescale, and $\delta$ is a nondimensional parameter for evaporation. This type of evaporative flux was first proposed by Burelbach et al. [6] in a one-sided evaporating thin film model where the dynamics of the fluid is assumed decoupled from the evolution of the vapor, and was later investigated by Ajaev [2, 1] for a more detailed evaporation model. Later a revised version of (2.7) was further studied in [17]. More recently the influence of surfactant was also included in the evaporative mass flux in [24].

In 2014, Peng et al. [20] derived a tear film break-up model with instabilities driven by evaporation effects by treating the lipid layer on the fluid tear film as a barrier to local evaporation from the underlying tear film. Assuming that the evolution of the lipid layer is static compared to the dynamics of the aqueous film, they imposed a fixed spatially varying profile to approximate the local variations in the lipid concentration. Since reduced amounts of the lipid cause locally elevated liquid evaporation rates relative to that of the film in surrounding regions [19, 5], they included an additional mass-transfer resistance term to account for the counteraction due to the lipid concentration. The resulting evaporative mass flux $J_e$ is then obtained by solving a coupled system for $J_e(x,t)$ and the temperature of the liquid-vapor interface. In order to simplify the evaporative flux term, we consider the same influence of lipid layer as an obstruction to evaporation but ignore the latent heat of water vaporization, and write the evaporative flux as

$$J_e = \frac{P_s e^{\delta p} - P_\infty}{R_L(x) + R_G} \sim \frac{E}{R_L(x) + R_G} \quad \text{for } \delta \to 0 \text{ with } E = P_s - P_\infty, \tag{2.8}$$

where $R_L(x)$ measures the mass-transfer resistance through the tear film lipid layer, which depends on the lipid concentration, $R_G$ represents the mass-transfer resistance in ambient air, $P_s$ and $P_\infty$ are nondimensional saturation vapor pressure at cornea and in the environment. This form of evaporative flux is comparable to (2.7) when the contribution from dynamical pressure $p$ to the evaporation effect is negligible and the film thickness $h \ll \bar{K}$, and water evaporation from the tear film is obstructed by both resistances through the spatially dependent lipid layer and through the air phase. Consequently, using (2.6) and (2.8), the total non-conservative flux from evaporative and osmolarity weeping flows can be written conveniently in terms of an $\bar{S}(x)$ function as

$$-J_e + J_w = s - \left(\frac{E}{R_L(x) + R_G} + S_0\right) = s - \bar{S}(x). \tag{2.9}$$

To summarize, from (2.1) and (2.4) the nondimensional governing equations for the evolution of tear film thickness $h$ and the osmolarity $s$ in human eyes can be represented by

$$h_t = -\left[h^3 \left(h_{xx} + \frac{\epsilon}{h^3}\right)_x\right]_x - (\bar{S} - s), \tag{2.10a}$$

$$s_t = s_{xx} + \left(\frac{h_x}{h} - h^2 \left(h_{xx} + \frac{\epsilon}{h^3}\right)_x\right) s_x + \frac{s}{h}(\bar{S} - s), \tag{2.10b}$$

with the associated boundary conditions (1.2), where (2.10b) is obtained from applying the product rule to the time derivative term and substituting equation (2.10a) into equation (2.4).

It is interesting to note that both the fourth-order term due to surface tension and the second-order term due to conjoining pressure in (2.10a) are stabilizing. Since our major interest is to examine the mechanism of possible breakdowns of the tear film model (2.10), in this work we neglect the conjoining pressure by setting $\epsilon = 0$ in (2.10) and focus on the competition between the fourth-order regularizing term and the non-conservative contributions. To get a better understanding of the PDEs, we consider the generalized model (1.1) to explore the key features of tear-film break-up with power-law mobility functions for both conservative and non-conservative contributions. In particular, in order to regularize the non-conservative effects in the PDEs, modified versions of the non-conservative terms are considered with the original terms $(\bar{S} - s)$ in equations (2.10) multiplied by a regularizing factor $h^m$. It is worth noting that the physical model (2.10) for tear films with $\epsilon = 0$ corresponds to the case $(m, n) = (0, 3)$ of the generalized PDE system (1.1).

In the present work, the PDE system (1.1) is investigated from the perspective of both analytical and numerical studies. Specifically, it will be shown numerically in section 4 that the tear film model (2.10) exhibits a novel finite-time rupture-shock phenomenon, that is, at a critical point $x = x_c$, the film thickness $h(x_c, t) \to 0$, along with $|s_x(x_c)| \to \infty$ as a critical time $t_c$ is approached. The finite-time singularity phenomenon admitted in this model inspires us to investigate the modified PDE system (1.1) in other parameter ranges.

## 3. Regularity of solutions to the generalized model (1.1)

In this section, we shall show the regularity of solutions to (1.1) with proper initial data and that the existence of strong solutions depends on different ranges of parameters $m$ and $n$. The main result for global strong solutions to equations (1.1) with $(m, n)$ that satisfy $n = m + 1$, $3 \leq m < 4$ is stated in Theorem 1, and the local existence of strong solutions to (1.1) with $m, n \geq 0$ is presented in Theorem 2.

*3.1. Global strong solution for $3 \leq m < 4$ and $n = m + 1$*

**Theorem 1.** *Let $\bar{S}(x) \in L^\infty([0, L])$. Suppose that $n = m + 1$, $3 \leq m < 4$. For any integer $k \geq 2$, positive constants $\eta, \lambda > 0$, and initial data $h_0(x) \in H^k([0, L])$, $s_0(x) \in H^{k-2}([0, L])$ satisfying that $0 < \eta \leq h_0(x)$, $0 < \lambda \leq s_0(x) \leq \|\bar{S}(x)\|_\infty$. Then for any $T > 0$, there exist $h(x, t) \in L^\infty([0, T]; H^k([0, L])) \cap L^2([0, T]; H^{k+1}([0, L]))$ and $s(x, t) \in L^\infty([0, T]; H^{k-2}([0, L])) \cap L^2([0, T]; H^{k-1}([0, L]))$ being the strong solution of (1.1) with initial data $h_0, s_0$ and boundary condition (1.2). Moreover, $s$ satisfies*

$$0 < \lambda \leq s(x, t) \leq \|\bar{S}\|_\infty, \quad t \in [0, T]. \tag{3.1}$$

*There exist positive constants $h_m(T)$, $H_m(T)$ such that*

$$0 < h_m(T) \leq h(x, t) \leq H_m(T), \text{ for all } t \in [0, T], \tag{3.2}$$

*where $h_m(T)$, $H_m(T)$ depend only on $\eta$, $\lambda$, $\|h_0 s_0\|_1$, $\|h_0\|_1$, $\|\bar{S}\|_\infty$, $T$ and $H_m(T)$ is the solution of*

$$H = \sqrt{C_\eta + C_s T + C(h_0, s_0, \lambda) H^\tau} + C(h_0, s_0, \lambda) + 1. \tag{3.3}$$

*with $2(m - 3) < \tau < 2$. $C_\eta, C_s$, are constants depending separately on $\eta$, $\|\bar{S}\|_\infty$, and $C(h_0, s_0, \lambda)$ is a constant depending on $\lambda$, $\|h_0 s_0\|_1$, $\|h_0\|_1$.* ◇

The uniform lower and upper bound estimates (3.2) are crucial for the higher order estimates, which ensure the existence of the global strong solution. We will use some basic estimates and two a-priori assumptions to obtain the uniform lower and upper bound estimates (3.2) and then verify the a-priori assumptions. We will show that the condition $m < 4$ comes from the upper bound $H_m(T)$ and that the condition $3 \leq m$ comes from the lower bound $h_m(T)$. Then standard compactness arguments give the existence result for global strong solutions.

If we consider strong solutions existing for local time, (3.2) can be easily obtained, which is important for the higher order estimate. Hence the conditions $3 \leq m < 4$, $n = m + 1$ can be broadened to $m \geq 0$, $n \geq 0$ and we state the existence result for local strong solution to the original model with parameters $m \geq 0$, $n \geq 0$ as a byproduct in Theorem 2.

For the case $n = m + 1$, the original model (1.1) becomes

$$h_t = -h^m(\bar{S} - s) - (h^{m+1} h_{xxx})_x, \tag{3.4a}$$

$$s_t = s_{xx} + \left(\frac{h_x}{h} - h^m h_{xxx}\right) s_x + s(\bar{S} - s) h^{m-1}, \tag{3.4b}$$

with a single system parameter $m$.

Now we start to prove the main result Theorem 1. The key point is to obtain the upper bound and positive lower bound of $h$, which is shown in Step 1 and Step 2 separately. For the following analysis, we denote $\|\cdot\|_p$ as the standard norm for $L^p([0,L])$, $1 \leq p \leq \infty$.

*Proof of Theorem 1.* The strategy of the proof is to first get some a-priori estimates under the assumption $h \geq 0$ in Step 1. We also obtain the upper bound $H_m(T)$ in Step 1, which requires $m < 4$. Then we verify the a-priori assumption by obtaining a positive lower-bound $h_m(T)$ of $h$ in Step 2, which requires $m \geq 3$. In Step 3, we can use the lower-bound $h_m(T)$ to obtain higher order a-priori estimates.

Step 1. Basic a-priori estimates. For any $T > 0$, we use the assumption

$$h(x,t) \geq 0, \text{ for all } t \in [0,T]. \tag{3.5}$$

Notice that $h \geq 0$ by (3.5) and that $\lambda \leq s_0(x) \leq \|\bar{S}(x)\|_\infty$. When it is the first time $s = \lambda$ or $s = \|\bar{S}(x)\|_\infty$, we have $s_x = 0$. For the last term on the righthandside of (3.4b), we have

$$s(\bar{S} - s)h^{m-1}|_{s=\lambda} = \lambda(\bar{S} - \lambda)h^{m-1} \geq 0,$$

$$s(\bar{S} - s)h^{m-1}|_{s=\|\bar{S}(x)\|_\infty} = \|\bar{S}(x)\|_\infty(\bar{S} - \|\bar{S}(x)\|_\infty)h^{m-1} \leq 0,$$

Then from (3.4b) and the weak form of maximal principle, $\lambda \leq s_0(x) \leq \|\bar{S}(x)\|_\infty$ implies

$$\lambda \leq s(x,t) \leq \|\bar{S}(x)\|_\infty, \tag{3.6}$$

which gives (3.1). This estimate together with (1.3) and (1.4) shows that

$$\int_0^L h \, dx \leq C(\|h_0 s_0\|_1, \lambda), \text{ for any } t \in [0,T], \tag{3.7}$$

and

$$\int_0^T \int_0^L h^m(\bar{S} - s) \, dx \, dt \leq C(\|h_0\|_1). \tag{3.8}$$

Moreover, multiplying (3.4a) by $\frac{1-m}{h^m}$ and integrating by parts lead to

$$\frac{d}{dt} \int_0^L \frac{1}{h^{m-1}} \, dx = (m-1) \int_0^L (\bar{S} - s) \, dx - m(m-1) \int_0^L h_{xx}^2 \, dx. \tag{3.9}$$

Denote $C_\eta = \frac{1}{\eta^{m-1}}$. From (3.6) and (3.9), we know

$$\int_0^T \int_0^L h_{xx}^2 \, dx \, dt \leq C_\eta + C_s T, \tag{3.10}$$

and

$$\int_0^L \frac{1}{h^{m-1}} \, dx \leq C_\eta + C_s T, \text{ for any } t \in [0,T], \tag{3.11}$$

where $C_s$ is a constant depending only on $\|\bar{S}\|_\infty$.

Furthermore, we turn to estimate $\|h\|_\infty$ and $\|h_x\|_{L^\infty([0,T];L^2([0,L]))}$. We need another a-priori assumption

$$\|h\|_\infty \leq H, \tag{3.12}$$

where $H$ will be determined later. Multiplying (3.4a) by $-h_{xx}$ and integrating by parts lead to

$$\frac{d}{dt} \int_0^L \frac{1}{2} h_x^2 \, dx = \int_0^L h^m h_{xx}(\bar{S} - s) \, dx - \int_0^L h_{xxx}^2 h^{m+1} \, dx.$$

Thus by Young's inequality and (3.10), we have

$$\int_0^L \frac{1}{2} h_x^2 \, dx \leq \int_0^T \int_0^L h^m h_{xx}(\bar{S} - s) \, dx \, dt, \tag{3.13}$$

$$\leq \int_0^T \int_0^L h_{xx}^2 \, dx \, dt + \int_0^T \int_0^L h^{2m}(\bar{S} - s)^2 \, dx \, dt,$$

$$\leq C_\eta + C_s T + I, \text{ for any } t \in [0, T],$$

where $I := \int_0^T \int_0^L h^{2m}(\bar{S} - s)^2 \, dx \, dt$ needs to be treated in detail.

Let $m, n, \gamma \geq 1$ satisfying

$$\frac{1}{m} + \frac{1}{n} + \frac{1}{\gamma} = 1. \tag{3.14}$$

We will choose the values of $m, n, \gamma$ in the following. Then by Hölder's inequality and (3.7), we have

$$\int_0^L h^{2m}(\bar{S} - s)^2 \, dx = \int_0^1 [h^m(\bar{S} - s)]^{\frac{1}{m}} h^{2m - \frac{m}{m} - \frac{1}{n}}(\bar{S} - s)^{2 - \frac{1}{m}} h^{\frac{1}{n}} \, dx$$

$$\leq C(\|h_0 s_0\|_1, \lambda) \left( \int_0^L h^m(\bar{S} - s) \, dx \right)^{\frac{1}{m}} \left( \int_0^L h^{(2m - \frac{m}{m} - \frac{1}{n})\gamma} \, dx \right)^{\frac{1}{\gamma}}.$$

This, together with Hölder's inequality and (3.8), shows that

$$I \leq C(\|h_0 s_0\|_1, \lambda) \left[ \int_0^T \int_0^L h^m(\bar{S} - s) \, dx \, dt \right]^{\frac{1}{m}} \left[ \int_0^T \left( \int_0^L h^{(2m - \frac{m}{m} - \frac{1}{n})\gamma} \, dx \right)^{\frac{m'}{\gamma}} dt \right]^{\frac{1}{m'}}$$

$$\leq C(h_0, s_0, \lambda) \left[ \int_0^T \left( \int_0^L h^{(2m - \frac{m}{m} - \frac{1}{n})\gamma} \, dx \right)^{\frac{m'}{\gamma}} dt \right]^{\frac{1}{m'}} \tag{3.15}$$

where $m' \geq 1$ satisfying

$$\frac{1}{m} + \frac{1}{m'} = 1. \tag{3.16}$$

Here and in the following, we denote $C(h_0, s_0, \lambda)$ as a constant depending only on $\|h_0 s_0\|_1, \|h_0\|_1, \lambda$.

For any constant $0 \leq \tau < 2$. Noticing the a-priori assumption (3.12) and (3.15), we know

$$I \leq H^\tau C(\|h_0 s_0\|_1, \lambda) \left[ \int_0^T \left( \int_0^L h^{(2m - \frac{m}{m} - \frac{1}{n} - \tau)\gamma} \, dx \right)^{\frac{m'}{\gamma}} dt \right]^{\frac{1}{m'}}. \tag{3.17}$$

On one hand, from the Gagliardo-Nirenberg interpolation inequality, we know

$$\|h\|_{(2m - \frac{m}{m} - \frac{1}{n} - \tau)\gamma} \leq c \|h\|_1^{1-\theta} \|h_{xx}\|_2^\theta + c \|h\|_1,$$

$$\leq C(h_0, s_0, \lambda) \|h_{xx}\|_2^\theta + C(h_0, s_0, \lambda), \text{ for any } t \in [0, T], \tag{3.18}$$

where the index $\theta$ satisfies

$$\theta = \frac{2}{5} \left[ 1 - \frac{1}{(2m - \frac{m}{m} - \frac{1}{n} - \tau)\gamma} \right]. \tag{3.19}$$

On the other hand, from the relations of (3.14) and (3.16), we know $0 \leq \frac{m'}{\gamma} \leq 1$. Therefore, from (3.17) and (3.18), we have

$$I \leq H^\tau C(h_0, s_0, \lambda) \left[ \left( \int_0^T \|h_{xx}\|_2^{\theta m'(2m - \frac{m}{m} - \frac{1}{n} - \tau)} dt \right)^{\frac{1}{m'}} + 1 \right] \tag{3.20}$$

Since we have uniform bound (3.10), it remains to show that

$$\theta m' \left( 2m - \frac{m}{m} - \frac{1}{n} - \tau \right) < 2. \tag{3.21}$$

Using relations (3.14), (3.16) and (3.19), this reduces to

$$\frac{m}{m-1}\Big[(2m - \frac{m}{m} - \frac{1}{n} - \tau) - (1 - \frac{1}{m} - \frac{1}{n})\Big] < 5, \tag{3.22}$$

which is

$$\begin{cases} \text{case 1: } \frac{6-m}{6-2m+\tau} < m,\ 6 - 2m + \tau > 0; \text{ or} \\ \text{case 2: } m < \frac{6-m}{6-2m+\tau},\ 6 - 2m + \tau < 0; \text{ or} \\ \text{case 3: } 6 - m < 0,\ 6 - 2m + \tau = 0. \end{cases} \tag{3.23}$$

Notice $m > 1$ and $0 \leq \tau < 2$. Only case 1 can happen, which becomes the minimum requirement to guarantee (3.21)

$$m < 3 + \frac{\tau}{2}. \tag{3.24}$$

Thus (3.20) becomes

$$I \leq C(h_0, s_0, \lambda) H^\tau. \tag{3.25}$$

This, together with (3.13), gives

$$\int_0^L h_x^2 \, \mathrm{d}x \leq C_\eta + C_s T + C(h_0, s_0, \lambda) H^\tau, \text{ for any } t \in [0, T], \tag{3.26}$$

from which, we also know $h$ is continuous. Assume $h$ achieves its minimal value, denoted as $h_{\min}$, at $x_c$. Since (3.7) shows $h_{\min} \leq C(h_0, s_0, \lambda)$, we have, from (3.26),

$$\|h\|_\infty \leq |\int_{x_c}^x h_x(s)\,\mathrm{d}s| + h_{\min} \leq \sqrt{C_\eta + C_s T + C(h_0, s_0, \lambda) H^\tau} + C(h_0, s_0, \lambda), \text{ for any } t \in [0, T]. \tag{3.27}$$

Finally, we verify the a-priori assumption (3.12). Let us choose $H$ being the solution of

$$H = \sqrt{C_\eta + C_s T + C(h_0, s_0, \lambda) H^\tau} + C(h_0, s_0, \lambda) + 1. \tag{3.28}$$

In fact, since $\tau < 2$, we can always find a solution, denoted as $H_m(T)$, to this equation, which depends only on $C_\eta$, $C(h_0, s_0, \lambda)$, $C_s$ and $T$. Then from (3.27), we have

$$\|h\|_\infty \leq \sqrt{C_\eta + C_s T + C(h_0, s_0, \lambda) H^\tau} + C(h_0, s_0, \lambda) < H_m(T), \text{ for any } t \in [0, T],$$

which verifies the a-priori assumption (3.12). Besides, for any $m < 4$, there exists $\delta > 0$ such that $m < 4 - \delta$. Then we can choose $2 - 2\delta < \tau < 2$, which implies $m < 4 - \delta < 3 + \frac{\tau}{2}$. Therefore, we obtain

$$\|h_x\|_{L^\infty([0,T];L^2([0,L]))} \leq C_{h_0, s_0, T}, \tag{3.29}$$

and

$$\|h\|_{L^\infty([0,T]\times[0,L])} \leq H_m(T), \tag{3.30}$$

where $C_{h_0, s_0, T}$ is a constant depending only on $\eta$, $\|h_0\|_1$, $\|\bar{S}\|_\infty$ and $T$.

Step 2. Positive lower bound for $h$. First from (3.29), (3.30) and the Sobolev embedding $H^1([0, L]) \hookrightarrow C^{\frac{1}{2}}([0, L])$, we know

$$h(x) \leq h_{\min} + C_{h_0, s_0, T} |x - x_c|^{\frac{1}{2}}. \tag{3.31}$$

This, together with (3.11), shows that

$$\int_0^L \frac{1}{(h_{\min} + C_{h_0, s_0, T}|x - x_c|^{\frac{1}{2}})^{m-1}} \,\mathrm{d}x$$
$$\leq \int_0^L \frac{1}{h^{m-1}} \,\mathrm{d}x \leq C_\eta + C_s T, \qquad \text{for any } t \in [0, T]. \tag{3.32}$$

If $m = 3$, we have
$$\ln \frac{h_{\min}^2 + C_{h_0,s_0,T}}{h_{\min}^2} = \int_0^{\frac{L}{2}} \frac{1}{h_{\min}^2 + C_{h_0,s_0,T}\, x}\, \mathrm{d}x \leq 2C_\eta + 2C_s T,$$
where $C_{h_0,s_0,T}$ is a constant depending only on $\eta, \|h_0\|_1, \|\bar{S}\|_\infty$ and $T$. Hence
$$1 + \frac{C_{h_0,s_0,T}}{h_{\min}^2} \leq e^{2C_\eta + 2C_s T},$$
which gives
$$h_{\min} \geq \left(\frac{C_{h_0,s_0,T}}{e^{2C_\eta + 2C_s T} - 1}\right)^{\frac{1}{2}}. \tag{3.33}$$

If $m > 3$, we have
$$\frac{\varepsilon}{(h_{\min} + C_{h_0,s_0,T}\, \varepsilon^{\frac{1}{2}})^{m-1}} \leq \int_0^\varepsilon \frac{1}{(h_{\min} + C_{h_0,s_0,T}\, \varepsilon^{\frac{1}{2}})^{m-1}}\, \mathrm{d}x \leq 2C_\eta + 2C_s T,$$
which gives
$$h_{\min} \geq \left(\frac{\varepsilon}{2C_\eta + 2C_s T}\right)^{\frac{1}{m-1}} - C_{h_0,s_0,T}\, \varepsilon^{\frac{1}{2}}. \tag{3.34}$$

Since $\frac{1}{m-1} < \frac{1}{2}$, we can choose $\varepsilon$ small enough such that the $\left(\frac{\varepsilon}{2C_\eta + 2C_s T}\right)^{\frac{1}{m-1}} - C_{h_0,s_0,T}\, \varepsilon^{\frac{1}{2}} > 0$. Combining the above two cases, we know, for $m \geq 3$, there exist a positive $h_m(T) > 0$ depending only on $\eta, \|h_0\|_1, \|\bar{S}\|_\infty$ and $T$, such that
$$h(x,t) \geq h_m(T) > 0, \text{ for any } t \in [0,T],$$
which gives (3.2) and verifies the a-priori assumption (3.5).

Step 3. Higher order a-priori estimates. Now we can use (3.2) to obtain higher order estimate.

First, we try to estimate $\|h_{xx}\|_{L^\infty([0,T],L^2([0,L]))}$. Denote $h^{(k)} := \frac{\partial^k h}{\partial x^k}$. Multiply (3.4a) by $h^{(4)}$ and integrate by parts. From Young's inequality, we have

$$\frac{\mathrm{d}}{\mathrm{d}t} \int_0^L \frac{1}{2} h_{xx}^2\, \mathrm{d}x = \int_0^L -h^m(\bar{S} - s)h^{(4)} - h^{m+1}(h^{(4)})^2 - (h^{m+1})_x h_{xxx} h^{(4)}\, \mathrm{d}x \tag{3.35}$$
$$\leq \int_0^L -h_m^{m+1}(h^{(4)})^2 + \varepsilon(h^{(4)})^2 + c(\varepsilon) H_m^{2m} \|\bar{S}\|_\infty^2\, \mathrm{d}x + I_1,$$

where $I_1 := \int_0^L -(h^{m+1})_x h_{xxx} h^{(4)}\, \mathrm{d}x$. We want to use $\int_0^L -h_m^{m+1}(h^{(4)})^2\, \mathrm{d}x$ to control $I_1$.

$$I_1 = -\frac{1}{2}\int_0^L (h^{m+1})_x (h_{xxx}^2)_x\, \mathrm{d}x = \frac{1}{2}\int_0^L (h^{m+1})_{xx} h_{xxx}^2\, \mathrm{d}x$$
$$\leq C(H_m) \int_0^L h_{xx} h_{xxx}^2\, \mathrm{d}x + C(H_m)\int_0^L h_x^2 h_{xxx}^2\, \mathrm{d}x =: I_2 + I_3.$$

Next, keeping in mind we have uniform bound (3.29), we use Gagliardo-Nirenberg interpolation inequality to estimate $I_2, I_3$. For $I_2$, we obtain
$$\|h_{xx}\|_\infty \leq \|h_x\|_2^{1-\theta_2} \|h^{(4)}\|_2^{\theta_2} + C_{h_0,s_0,T}, \text{ for } \theta_2 = \frac{1}{2}, \tag{3.36}$$
and
$$\|h_{xxx}\|_2 \leq \|h_x\|_2^{1-\theta_1} \|h^{(4)}\|_2^{\theta_1} + C_{h_0,s_0,T}, \text{ for } \theta_1 = \frac{2}{3}. \tag{3.37}$$

Since $\theta_2 + 2\theta_1 = \frac{11}{6} < 2$, (3.36) and (3.37), together with Young's inequality, show that
$$I_2 \leq C(H_m)\|h_{xx}\|_\infty \|h_{xxx}\|_2^2 \leq \varepsilon \|h^{(4)}\|_2^2 + C_{h_0,s_0,T}, \tag{3.38}$$

where $C_{h_0,s_0,T}$ is a constant depending only on $\eta, \|h_0\|_1, \|\bar{S}\|_\infty$ and $T$. Similarly, for $I_3$, we obtain

$$\|h_x\|_\infty \leq \|h_x\|_2^{1-\theta_4}\|h^{(4)}\|_2^{\theta_4} + C_{h_0,s_0,T}, \text{ for } \theta_4 = \frac{1}{6}, \tag{3.39}$$

and

$$\|h_{xxx}\|_2 \leq \|h_x\|_2^{1-\theta_3}\|h^{(4)}\|_2^{\theta_3} + C_{h_0,s_0,T}, \text{ for } \theta_3 = \frac{2}{3}. \tag{3.40}$$

Since $2\theta_3 + 2\theta_4 = \frac{5}{3} < 2$, (3.39) and (3.40), together with Young's inequality, show that

$$I_3 \leq C(H_m)\|h_x\|_\infty^2 \|h_{xxx}\|_2^2 \leq \varepsilon\|h^{(4)}\|_2^2 + C_{h_0,s_0,T}, \tag{3.41}$$

where $C_{h_0,s_0,T}$ is a constant depending only on $\eta, \|h_0\|_1, \|\bar{S}\|_\infty$ and $T$. Combining (3.38), (3.41) with (3.35), we have

$$\frac{d}{dt}\int_0^L h_{xx}^2\,dx + \int_0^L h_m^{m+1}(h^{(4)})^2\,dx \leq C_{h_0,s_0,T}, \tag{3.42}$$

which gives that

$$\|h_{xx}\|_{L^\infty([0,T],L^2([0,L]))} \leq C_{h_0,s_0,T}, \tag{3.43}$$

$$\|h^{(4)}\|_{L^2([0,T],L^2([0,L]))} \leq C_{h_0,s_0,T}, \tag{3.44}$$

where $C_{h_0,s_0,T}$ is a constant depending only on $\eta, \|h_0\|_1, \|\bar{S}\|_\infty$ and $T$.

Second, in order to get higher order $h$-estimates, we now need to obtain $s$-estimate. Multiplying (3.4b) by $s$ and integrate by parts, we have

$$\frac{d}{dt}\int_0^L \frac{1}{2}s^2\,dx \tag{3.45}$$

$$= -\int_0^L s_x^2\,dx + \int_0^L \left(\frac{h_x}{h} - h^m h_{xxx}\right)ss_x + s^2(\bar{S}-s)h^{m-1}\,dx$$

$$\leq -\int_0^L s_x^2\,dx + \varepsilon\int_0^L s_x^2\,dx + C(\varepsilon,\|\bar{S}\|_\infty)\int_0^L \left(\frac{h_x}{h} - h^m h_{xxx}\right)^2\,dx + \int_0^L s^2(\bar{S}-s)h^{m-1}\,dx$$

$$\leq -\int_0^L s_x^2\,dx + \varepsilon\int_0^L s_x^2\,dx + C(\varepsilon,\|\bar{S}\|_\infty)\int_0^L \left(2\frac{(h_x)^2}{h_m^2} + 2H_m^{2m}(h_{xxx})^2\right)\,dx + C_{h_0,s_0,T},$$

where we used Young's inequality and (3.2). Notice estimates (3.43) and (3.44). Integrating $t$ from 0 to $T$, (3.45) yields

$$\int_0^L s^2\,dx + \int_0^T\int_0^L s_x^2\,dx \leq C_{h_0,s_0,T} \tag{3.46}$$

for any $t \in [0,T]$, where $C_{h_0,s_0,T}$ is a constant depending only on $\eta, \|h_0\|_1, \|\bar{S}\|_\infty$ and $T$.

Third, we turn to estimate $\|h_{xxx}\|_{L^\infty([0,T],L^2([0,L]))}$. Multiply (3.4a) by $h^{(6)}$ and integrate by parts. We have

$$\frac{d}{dt}\int_0^L \frac{1}{2}h_{xxx}^2\,dx = \int_0^L h^m(\bar{S}-s)h^{(6)} - (h^{m+1}h_{xxx})_{xx}h^{(5)}\,dx \tag{3.47}$$

$$= \int_0^L -(h^m(\bar{S}-s))_x h^{(5)}\,dx - \int_0^L h^{m+1}(h^{(5)})^2\,dx$$

$$- \int_0^L 2(h^{m+1})_x h^{(4)}h^{(5)} + (h^{m+1})_{xx}h_{xxx}h^{(5)}\,dx,$$

$$\leq \varepsilon\int_0^L (h^{(5)})^2\,dx + C(\varepsilon)C_{h_0,s_0,T} - \int_0^L h_m^{m+1}(h^{(5)})^2\,dx$$

$$+ \int_0^L (h^{m+1})_{xx}(h^{(4)})^2\,dx - \int_0^L (h^{m+1})_{xx}h_{xxx}h^{(5)}\,dx,$$

where we used Young's inequality, (3.2), (3.29) and (3.46) in the last inequality. Denote $R := \int_0^L (h^{m+1})_{xx}(h^{(4)})^2 \, dx$ and $R_3 := -\int_0^L (h^{m+1})_{xx} h_{xxx} h^{(5)} \, dx$.

$$R = \int_0^L (h^{m+1})_{xx}(h^{(4)})^2 \, dx$$
$$\leq C(H_m) \int_0^L h_{xx}(h^{(4)})^2 \, dx + C(H_m) \int_0^L h_x^2 (h^{(4)})^2 \, dx$$
$$=: R_1 + R_2.$$

We now use the term $-\int_0^L h_m^{m+1}(h^{(5)})^2 \, dx$ to control $R_1$, $R_2$ and $R_3$.

For $R_1$, using Gagliardo-Nirenberg interpolation inequality and keeping in mind the uniform bound (3.43), we obtain
$$\|h_{xx}\|_\infty \leq \|h_{xx}\|_2^{1-\theta_2} \|h_{xxx}\|_2^{\theta_2} + C_{h_0,s_0,T}, \text{ for } \theta_2 = \frac{1}{2}, \tag{3.48}$$

and
$$\|h^{(5)}\|_2 \leq \|h_{xx}\|_2^{1-\theta_1} \|h^{(5)}\|_2^{\theta_1} + C_{h_0,s_0,T}, \text{ for } \theta_1 = \frac{2}{3}. \tag{3.49}$$

Then (3.48) and (3.49) show that
$$R_1 \leq C(H_m)\|h_{xx}\|_\infty \|h^{(4)}\|_2^2$$
$$\leq C_{h_0,s_0,T}\|h_{xxx}\|_2^{\theta_2} \|h^{(5)}\|_2^{2\theta_1} + C_{h_0,s_0,T}\|h_{xxx}\|_2^{\theta_2} + C_{h_0,s_0,T}\|h^{(5)}\|_2^{2\theta_1} + C_{h_0,s_0,T} \tag{3.50}$$

where $C_{h_0,s_0,T}$ is a constant depending only on $\eta, \|h_0\|_1, \|\bar{S}\|_\infty$ and $T$. From Young's inequality, we have
$$\|h_{xxx}\|_2^{\theta_2} \|h^{(5)}\|_2^{2\theta_1} \leq \varepsilon \|h^{(5)}\|_2^{2p\theta_1} + C(\varepsilon)\|h_{xxx}\|_2^{q\theta_2},$$

with $q = 3$, $p = \frac{3}{2}$. Since $q\theta_2, 2\theta_1 < 2$, from Young's inequality again, we obtain
$$R_1 \leq \varepsilon \|h^{(5)}\|_2^2 + C_{h_0,s_0,T}\|h_{xxx}\|_2^2 + C_{h_0,s_0,T}. \tag{3.51}$$

For $R_2$, using (3.43) and Sobolev interpolation inequality, we have
$$R_2 = C(H_m) \int_0^L h_x^2(h^{(4)})^2 \, dx \leq C(H_m)\|h_x\|_\infty^2 \|h^{(4)}\|_2^2 \tag{3.52}$$
$$\leq C(H_m)\|h_{xx}\|_2^2 \|h^{(4)}\|_2^2 \leq C_{h_0,s_0,T}\|h^{(4)}\|_2^2$$
$$\leq \varepsilon \|h^{(5)}\|_2^2 + C_{h_0,s_0,T}.$$

For $R_3$, from Young's inequality and Sobolev interpolation inequality, we have
$$R_3 = -\int_0^L (h^{m+1})_{xx} h_{xxx} h^{(5)} \, dx \tag{3.53}$$
$$\leq \varepsilon \|h^{(5)}\|_2^2 + \|h_{xxx}\|_\infty^2 \int_0^L (h^{m+1})_{xx}^2 \, dx$$
$$\leq \varepsilon \|h^{(5)}\|_2^2 + (\varepsilon \|h^{(5)}\|_2^2 + C_{h_0,s_0,T})(\|h_{xx}\|_2^2 + \|h_x\|_4^4)$$
$$\leq \varepsilon \|h^{(5)}\|_2^2 + C_{h_0,s_0,T}. \tag{3.54}$$

where we used (3.43).

Combining (3.51), (3.52) and (3.53) with (3.47), we obtain
$$\frac{d}{dt} \int_0^L h_{xxx}^2 \, dx + \int_0^L h_m^{m+1}(h^{(5)})^2 \, dx \leq C_{h_0,s_0,T}\|h_{xxx}\|_2^2 + C_{h_0,s_0,T}. \tag{3.55}$$

This, together with Grönwall inequality, implies
$$\|h_{xxx}\|_{L^\infty([0,T],L^2([0,L]))} \leq C_{h_0,s_0,T}, \tag{3.56}$$

$$\|h^{(5)}\|_{L^2([0,T],L^2([0,L]))} \leq C_{h_0,s_0,T}, \tag{3.57}$$

where $C_{h_0,s_0,T}$ is a constant depending only on $\eta, \|h_0\|_1, \|\bar{S}\|_\infty$ and $T$.

Finally, we turn to estimate $\|s_x\|_{L^\infty([0,T],L^2([0,L]))}$. Multiply (3.4b) by $-s_{xx}$ and integrate by parts. From Young's inequality, (3.6) and (3.2), we have

$$\frac{d}{dt}\int_0^L \frac{1}{2}s_x^2\,dx = \int_0^L -s_{xx}^2 - \left(\frac{h_x}{h} - h^m h_{xxx}\right)s_x s_{xx} - ss_{xx}(\bar{S}-s)h^{m-1}\,dx \tag{3.58}$$

$$\leq \int_0^L -s_{xx}^2\,dx + \varepsilon\int_0^L s_{xx}^2\,dx + \frac{1}{2}\int_0^L \left(\frac{h_x}{h} - h^m h_{xxx}\right)_x s_x^2\,dx + C_{h_0,s_0,T},$$

where we want to use the first term on the righthandside to control $P := \frac{1}{2}\int_0^L (\frac{h_x}{h} - h^m h_{xxx})_x s_x^2\,dx$. From (3.29), (3.43) and (3.2), we know

$$P \leq C_{h_0,s_0,T}\left[\int_0^L h_{xx}s_x^2 + h_x^2 s_x^2 + h_x h_{xxx} s_x^2 + h^{(4)} s_x^2\,dx\right] \tag{3.59}$$

$$\leq C_{h_0,s_0,T}\left[(\|h_{xx}\|_\infty + \|h_x\|_\infty^2)\int_0^L s_x^2\,dx + \|h_x\|_\infty\int_0^L h_{xxx}s_x^2\,dx + \int_0^L h^{(4)}s_x^2\,dx\right]$$

$$\leq \varepsilon\int_0^L s_{xx}^2\,dx + C_{h_0,s_0,T}\left(\int_0^L h_{xxx}s_x^2\,dx + \int_0^L h^{(4)}s_x^2\,dx\right) + C_{h_0,s_0,T}.$$

Now we estimate $\int_0^L h^{(4)}s_x^2\,dx$, the other term $\int_0^L h_{xxx}s_x^2\,dx$ is same.

$$\int_0^L h^{(4)}s_x^2\,dx \leq C(\varepsilon)\int_0^L (h^{(4)})^q\,dx + \varepsilon\int_0^L s_x^{2p}, \tag{3.60}$$

with $p = \frac{5}{3}$, and $q = \frac{5}{2}$. Using Gagliardo-Nirenberg interpolation inequality and (3.46), we know

$$\|s_x\|_{2p} \leq \|s_{xx}\|_2^{\theta_1}\|s\|_2^{1-\theta_1} + C_{h_0,s_0,T},\ \text{for}\ \theta_1 = \frac{3}{5}, \tag{3.61}$$

and

$$\|h^{(4)}\|_q \leq \|h^{(5)}\|_2^{\theta_2}\|h_{xxx}\|_2^{1-\theta_2} + C_{h_0,s_0,T},\ \text{for}\ \theta_2 = \frac{11}{20}. \tag{3.62}$$

Thus we have

$$\int_0^L s_x^{2p}\,dx \leq C_{h_0,s_0,T}\int_0^L s_{xx}^2\,dx + C_{h_0,s_0,T}, \tag{3.63}$$

and since $q\theta_2 = \frac{55}{40} < 2$, we know

$$\int_0^L (h^{(4)})^q\,dx \leq C_{h_0,s_0,T}\int_0^L (h^{(5)})^2\,dx + C_{h_0,s_0,T}. \tag{3.64}$$

Combining (3.61), (3.63), (3.64) with (3.59), we know

$$P \leq \varepsilon\int_0^L s_{xx}^2\,dx + C_{h_0,s_0,T}\int_0^L (h^{(5)})^2\,dx + C_{h_0,s_0,T}.$$

This, together with (3.58), gives

$$\frac{d}{dt}\int_0^L s_x^2\,dx + \int_0^L s_{xx}^2\,dx \leq C_{h_0,s_0,T}\int_0^L (h^{(5)})^2\,dx + C_{h_0,s_0,T}. \tag{3.65}$$

Noticing (3.57), integrate (3.65) from 0 to $T$, we obtain

$$\|s_x\|_{L^\infty([0,T],L^2([0,L]))} \leq C_{h_0,s_0,T}, \tag{3.66}$$

$$\|s_{xx}\|_{L^2([0,T],L^2([0,L]))} \leq C_{h_0,s_0,T}, \tag{3.67}$$

where $C_{h_0,s_0,T}$ is a constant depending only on $\eta, \|h_0\|_1, \|\bar{S}\|_\infty$ and $T$.

We can use the same techniques to obtain any $k$-th order estimates and by standard compactness arguments, we can obtain the existence result for global strong solution to (3.4a), (3.4b). We omit the details here and complete the proof of Theorem 1. □

*3.2. Local strong solution with $m \geq 0$, $n \geq 0$*

Notice the condition $3 \leq m < 4$ comes from the uniform lower and upper bound estimate (3.2), which is crucial for the higher order estimate. If we consider strong solutions existing for local time, (3.2) can be easily obtained. Next, we state the local strong solution for the original model with index $m \geq 0$, $n \geq 0$.

**Theorem 2.** *Let $\bar{S}(x) \in L^\infty([0,L])$. For any integer $k \geq 2$, positive constants $\lambda, \eta > 0$, assume $m \geq 0$, $n \geq 0$ and initial data $h_0(x) \in H^k([0,L])$, $s_0(x) \in H^{k-2}([0,L])$ satisfying that $h_0(x) \geq \eta > 0$, $0 < \lambda \leq s_0(x) \leq \|\bar{S}(x)\|_\infty$. Then there exists $T_m > 0$, depending only on $\eta, \lambda, \|h_0\|_1, \|h_0 s_0\|_1, \|\bar{S}\|_\infty$, such that $h(x,t) \in L^\infty([0,T_m]; H^k([0,L])) \cap L^2([0,T_m]; H^{k+1}([0,L]))$ and $s(x,t) \in L^\infty([0,T_m]; H^{k-2}([0,L])) \cap L^2([0,T_m]; H^{k-1}([0,L]))$ are the strong solution of (1.1) with initial data $h_0, s_0$ and boundary condition (1.2). Moreover,*

$$0 < \lambda \leq s(x,t) \leq \|\bar{S}\|_\infty, \quad t \in [0, T_m]. \tag{3.68}$$

*There exist positive constants $h_m(T_m)$, $H_m(T_m)$ such that*

$$0 < h_m(T_m) \leq h(x,t) \leq H_m(T_m), \text{ for all } t \in [0, T_m], \tag{3.69}$$

*where $h_m, H_m$ depend only on $\eta, \lambda, \|h_0\|_1, \|h_0 s_0\|_1, \|\bar{S}\|_\infty$ and $T_m$.* ◇

Since for local solution we can obtain (3.2) easily, the techniques for a-priori estimates are the same as the proof of Theorem 1 and thus we omit the details here.

## 4. Dynamics of model (1.1): numerical study

It was analytically shown in the previous section that for parameters $n = m + 1$ and $3 \leq m < 4$ the global existence of solutions to (1.1) is guaranteed. In order to verify this conclusion, we conduct a series of numerical simulations for the generalized equations (1.1) with different values of mobility exponents $(m,n)$ to investigate various long-time behaviors of the solutions. Furthermore, simulations reveal that the model has rich dynamics resulting from its coupled strong nonlinearity. For example, interesting finite-time singularities are observed with $(m,n) = (0,3)$ which corresponds to the tear film break up model (2.10). We will also discuss the significance of the effective salt capacity $\bar{S}(x)$ to the existence of equilibrium solutions and the dynamics of the model.

*4.1. Finite-time singularities*

While the motivating tear film break-up model (2.10) successfully captures the key components in evaporating tear films, the instabilities driven by the locally elevated evaporation rates can lead to a novel finite-time rupture-shock phenomenon. Fig. 2 shows a typical numerical simulation for the evolution of film height and osmolarity from initial condition $h_0 = s_0 = 1$ with $\bar{S}$ given by $\bar{S}(x) = 50 - 48.8 \tanh(20(|x-1|-0.1))$.

In the early stage of the dynamics, it is shown in Fig. 2 $(a,b)$ that the film thickness $h$ decreases with the osmolarity $s$ increasing since the locally elevated evaporation effects are large enough to overcome the curvature-driven and osmotic healing flows. In the later stage, Fig. 2 $(d)$ shows that $s$ locally exceeds the prescribed $\bar{S}$, and the local minimum of $h$ splits into a pair of secondary rupture structures (Fig. 2 $(c)$); at the same time the rupture in $h$ is smoothed by both weak diffusive and capillary forces. As thinning proceeds, the film thickness in the neighborhood of the critical position $x_c$ approaches zero which leads to degenerate diffusion for the local salt concentration from (2.4). This causes the osmolarity to form a singular shock in finite-time with $|s_x(x_c,t)| \to \infty$ for $t \to t_c$ (Fig. 2 $(f)$), with the development of tear film rupture as $h(x_c,t) \to 0$ (Fig. 2 $(e)$). It is clear from Fig. 2$(d)$ that the osmolarity $s$ is bounded by $\|\bar{S}\|_\infty$ throughout the simulation as is predicted by (3.1).

This type of rupture-shock dynamics is comparable to the double shock solutions studied in [13, 22] in a model for a thin viscous film with insoluble surfactant. That PDE system for film thickness $h$ and surfactant concentration $\Gamma$ allows shock solutions for which both $h$ and $\Gamma_x$ have a jump while $\Gamma$ is continuous. In particular, Jensen and

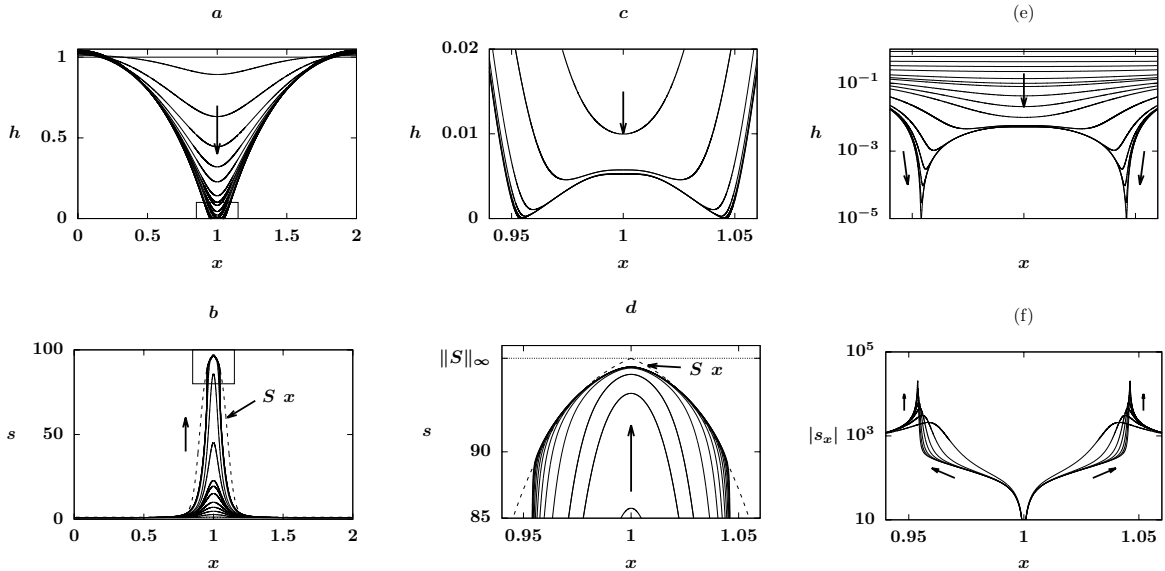

Figure 2: Evolution of $h$ and $s$ to equations (2.10) with $\epsilon = 0$, or equivalently, equations (1.1) with $(m, n) = (0, 3)$, starting from constant initial data $h_0 = s_0 = 1$ for $0 \leq x \leq 2$ driven by non-conservative flux with $\bar{S}(x) = 50 - 48.8\tanh(20(|x-1|-0.1))$ (plotted in dashed lines). Solutions profiles for $h$ and $s$ are shown in (a) and (b), with zoom-in plots in (c) and (d). In (e) and (f) solutions $h$ and $|s_x|$ are plotted on log scale, showing finite time rupture-shock singularity occurring at $x_c$ with $h(x_c) \to 0$ and $|s_x(x_c)| \to \infty$ as $t \to t_c$.

Grotberg [13] showed that severe film thinning behind the shock due to van der Waals can lead to film rupture. Different double shock and singular shock solutions for film thickness and particle volume fraction have also been investigated by Cook and Bertozzi for particle-laden thin films [8].

The nonlinear PDE system (1.1) is solved numerically using a fully implicit second-order finite difference method with an adaptive non-uniform grid. Specifically, we used the midpoint Keller-box method [14] to express the PDEs as a discrete system of first-order equations

$$h_t = -(h^n q)_x - h^m(\bar{S} - s) \qquad (hs)_t = (hw - h^n qs)_x, \tag{4.1a}$$

with

$$w = s_x, \quad k = h_x, \quad p = k_x, \quad q = p_x, \tag{4.1b}$$

where the second equation maintains the conservation of local salt mass from equation (2.4). To capture the finite-time rupture in $h$ and shock in $s$ that occur simultaneously with high resolutions, we used a classical moving mesh algorithm with a tailored monitor function together with adaptive time-stepping to adaptively assign a high distribution of grid points near the singularity points. For more discussion and applications of moving mesh methods, we refer to [7].

The presence of finite-time singularities indicates that the tear film model (2.10) is problematic since the solution $(h, s)$ cannot be continued past the time of the first singularity. Our regularization of the non-conservative contributions by introducing the mobility parameter $m$ to the generalized model (1.1) is inspired by this observation. For the following simulations we keep $n = m + 1$ so that the parameters are consistent with those used in Theorem 1. We also apply the initial conditions $h_0 = s_0 = 1$ with domain size $L = 2$, and define the effective salt capacity $\bar{S}(x)$ as a step function with a shift coefficient $\xi > 0$ that defines the width of the elevated-evaporation-rate region,

$$\bar{S}(x) = \begin{cases} 100 & \text{for} \quad L/2 - \xi < x < L/2 + \xi, \\ 2 & \text{otherwise.} \end{cases} \tag{4.2}$$

This choice of $\bar{S}$ and initial data satisfies the requirement of the global existence theorem $s_0 \leq \|\bar{S}\|_\infty$ and provides a typical characterization of the dynamics in the model (2.10).

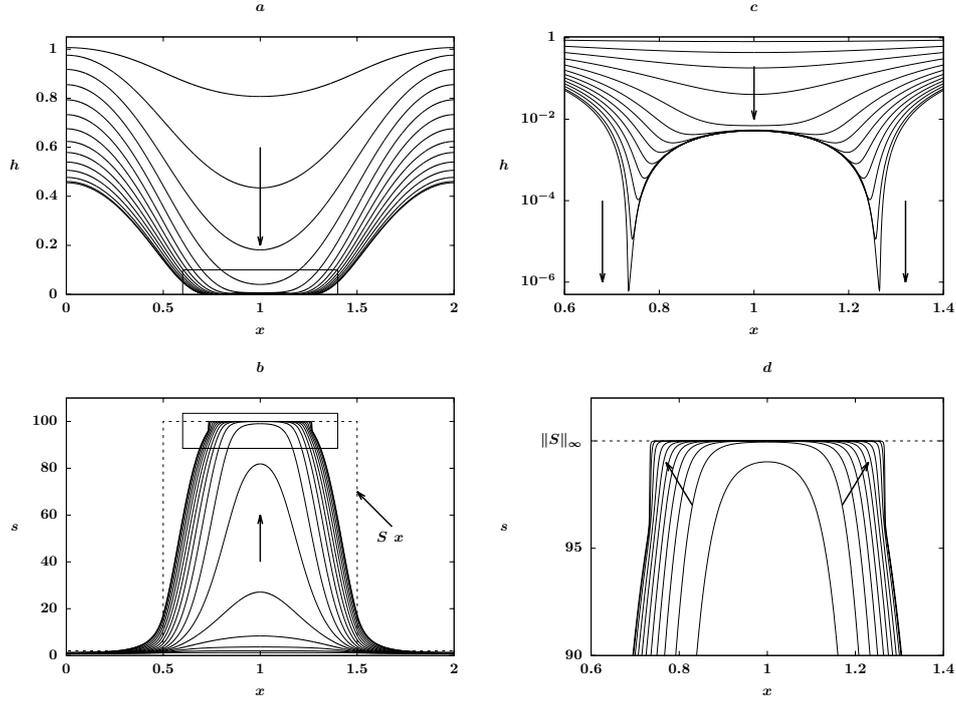

Figure 3: Numerical simulation of (1.1) with $(m,n) = (0.5, 1.5)$ and identical initial data and $\bar{S}$ (in dashed lines), (4.2) with $\xi = 0.5$, as in Fig. 1. Evolution of $h$ and $s$ are plotted in (a) and (b), with zoom-in plots in (c) and (d) showing that rupture-shock singularity occurs at a pair of points $x_c$ away from $x = 1$ with $h(x_c) \to 0$ and $|s_x(x_c)| \to \infty$ as $t \to t_c$, where $t_c \approx 0.063$.

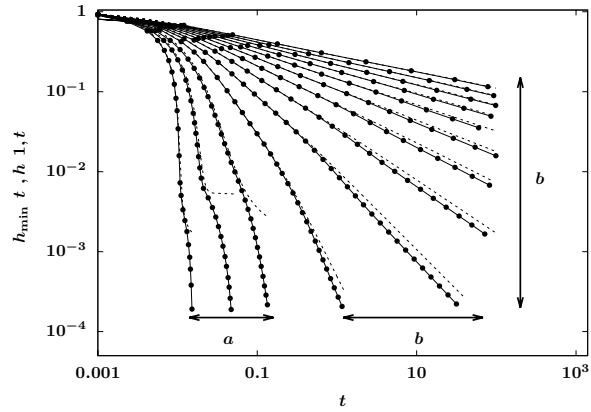

Figure 4: Plot of decreasing minimum film thickness $h_{\min}$ (plotted in dots) and $h(1,t)$ (plotted in dashed lines) starting from identical initial condition $h_0 = s_0 = 1$ with $n = m + 1$, $\bar{S}$ from (4.2) with $\xi = 0.5$ and over a range of $m$ values, $0 \leq m \leq 6$. Regime (a): $m = 0, 0.5, 1$; Regime (b): $m = 1.5, \cdots, 6$. In regime (a) finite-time rupture in $h$ develops, while infinite-time thinning occurs in regime (b). The deviation of $h_{\min}(t)$ from $h(1,t)$ in the later stage indicates that the minimum of the film thickness $h(x,t)$ is attained away from the center of the domain.

The PDE simulations shown in Fig. 3 with $(m, n) = (0.5, 1.5)$ and $\xi = 0.5$ suggest that weak regularization is not sufficient to prevent the finite-time singularities from happening. Similar to the dynamics presented in Fig. 2, rupture-shock phenomenon occurs in the later stage with $h(x_c) \to 0$ and $|s_x(x_c)| \to \infty$ as $t \to t_c \approx 0.063$. Since the width of the high capacity region in $\bar{S}$ increases from approximately 0.2 in Fig. 2 to $2\xi = 1$ in Fig. 3, the rupture hole in film thickness $h$ and the hyperosmotic region in $s$ are larger, while the secondary rupture-shock phenomenon that occurs away from the center of the domain is similar to the Fig. 2. Again the comparison between the profiles for $\bar{S}$ and the salt concentration $s$ in Fig. 3(d) emphasizes that $s$ does not exceed $\|\bar{S}\|_\infty = 100$ during the dynamical evolution.

Note that the distinct numerical simulations presented in Fig. 1 and Fig. 3 differ only in their choices of mobility exponents $(m, n)$, with all the other system parameters including the $\bar{S}$ profile being identical. Inspired by this observation, we explore the dynamics of (1.1) with $\bar{S}$ given by (4.2) with $\xi = 0.5$ and investigate the influences of various mobility exponents $m$ with $n = m + 1$. The time evolution of a sequence of PDE simulations with identical initial data is plotted in Fig. 4 with two different regimes. For $0 \le m \le 1$, localized finite-time rupture occurs at a point away from the origin similar to the dynamics shown in Fig. 3, while for $m > 1$ infinite-time non-uniform thinning is observed with the minimum film thickness $h_{\min} \to 0$ as $t \to \infty$, which is similar to the dynamics shown in Fig. 1. Specifically, the film thicknesses at the center of the domain $x = 1$ are plotted in comparison to the minimum film thickness. For regime $(a)$, as the finite-time singularity develops, $h_{\min}(t)$ quickly deviates from $h(1, t)$ and the difference between the two quantities grows exponentially as the critical time is approached, indicating the formation of the secondary singularities similar to the case shown in Fig. 3. These results support the conclusion drawn in Theorem 1 that strong solutions to (1.1) exist globally for $3 \le m < 4$ and $n = m + 1$. Moreover, the numerical result in Fig. 3 suggests that strong solutions to the model (1.1) exist until the first singularity occurs, which agrees with the local strong solution result in Theorem 2.

*4.2. Convergence to equilibrium and infinite time thinning*

It is shown in Fig. 4 that with $\bar{S}(x)$ from (4.2) and $\xi = 0.5$ one can separate the finite-time singularity regime of the solution behaviors from infinite-time thinning regime with various $(m, n)$ values. We shall then further investigate the long time behavior of the solutions of (1.1). In addition to the infinite time thinning, typical long-time behaviors of solutions $(h, s)$ of PDE system (1.1) may also include convergence to equilibrium solutions.

There is a possible equilibrium balance between the regularized non-conservative effects and the surface tension contributions in the PDE system. By setting the time-derivative terms in (1.1) equal to zero and applying the conservation of total mass of salt (2.5), we note that an equilibrium of the PDE system (1.1), $h_{\text{eq}}(x)$ and $s_{\text{eq}}(x)$, with initial data $(h_0, s_0)$ satisfies the differential equation system

$$\frac{d}{dx}\left(h_{\text{eq}}^n \frac{d^3 h_{\text{eq}}}{dx^3}\right) + h_{\text{eq}}^m(\bar{S} - s_{\text{eq}}) = 0, \tag{4.3a}$$

$$\frac{d^2 s_{\text{eq}}}{dx^2} + \left(\frac{1}{h_{\text{eq}}}\frac{dh_{\text{eq}}}{dx} - h_{\text{eq}}^{n-1}\frac{d^3 h_{\text{eq}}}{dx^3}\right)\frac{ds_{\text{eq}}}{dx} + s_{\text{eq}}(\bar{S} - s_{\text{eq}})h_{\text{eq}}^{m-1} = 0, \tag{4.3b}$$

$$\tag{4.3c}$$

subject to the constraint

$$\int_0^L h_{\text{eq}} s_{\text{eq}} \, dx = \int_0^L h_0 s_0 \, dx = Q_0 \tag{4.3d}$$

The existence of such equilibrium solutions depends on the profile of effective salt capacity $\bar{S}$ and other parameters. For instance, with initial condition constraint $Q_0 = 2$ and the form of $\bar{S}$ given by (4.2) with varying shift coefficient $\xi$, the equilibrium solutions to (1.1) with $(m, n) = (3.5, 4.5)$ are calculated via a continuation method and are plotted in Fig. 5. Note that for smaller $\xi$ in (4.2) the steady state $h_{\text{eq}}$ has a positive lower bound, while the minimum of $h_{eq}$ approaches zero at $x = 1$ with $\xi \sim 0.314$. This result suggests that steady states for (1.1) with $Q_0 = 2$ and $(m, n) = (3.5, 4.5)$ do not exist for $\xi > 0.314$ for $\bar{S}$ given by (4.2). It is interesting that while $\min h_{\text{eq}}$ is monotonically decreasing in terms of $\xi$, the profile of the equilibrium osmolarity $s_{\text{eq}}$ changes dramatically and $\max s_{\text{eq}}$ is not monotone in $\xi$, as is shown in the Fig. 5 (right).

Figure 6 depicts a typical simulation of the model (1.1) with the total mass of salt $Q_0 = 2$, $\xi = 0.2$ and $(m, n) = (3.5, 4.5)$, showing convergence of PDE solutions (dashed lines) to equilibrium solutions $(h_{\text{eq}}, s_{\text{eq}})$ (solid lines) in the long time. Moreover, numerical observations for $m > 1$ indicate that PDE solutions of (1.1) converge to the corresponding equilibrium, if it exists, which satisfies (4.3). If the equilibrium does not exist, we expect infinite-time non-uniform thinning with $h \to 0$ at a critical point $x_c$ as $t \to \infty$. With the shift coefficient $\xi = 0.3$, a sequence of PDE simulations starting from identical initial data are plotted in Figure 7 with $h_{\min}$ decreasing in time. Three distinct regimes are developed in this case: for region (a) with $0 \leq m \leq 1$ finite-time singularity develops similar to the case shown in Fig. 3; For region (c) with $m \geq 3.5$ the PDE solution converges to an equilibrium solution $(h_{eq}, s_{eq})$ similar to the dynamics in Fig. 6; While for region (b) with $1 < m \leq 3$ the minimum film thickness approaches zero as $t \to \infty$. Specifically, in the neighborhood of $x_c$ where $h \to 0$, the film thickness profile forms a nearly flat plateau with the corresponding $s \ll \bar{S}$. Therefore from (1.1a) the minimum film thickness is asymptotically determined by

$$\frac{d}{dt} h_{\min} \sim -\eta h_{\min}^m, \quad \text{where} \quad \eta = \bar{S}(x_c),$$

which leads to an estimate of the rate of change of $h_{\min}$ in time,

$$h_{\min}(t) \sim (c + \eta(m-1)t)^{-\frac{1}{m-1}}, \tag{4.4}$$

where $c$ is a constant that depends on other system parameters and initial conditions. The comparison of the direct PDE simulations against the prediction $h_{\min}(t) = O(t^{-\frac{1}{m-1}})$ as $t \to \infty$ for regime (b) in Fig. 7 shows good agreement with this estimate in the final stage as $h_{\min} \to 0$. The estimate in (4.4) also suggests that infinite-time thinning cannot happen for $m < 1$.

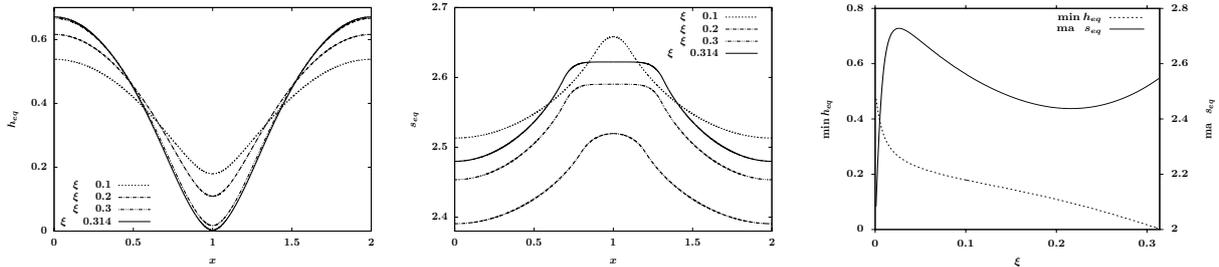

Figure 5: Equilibrium solutions $(h_{\text{eq}}, s_{\text{eq}})$ of (1.1) satisfying (4.3) with $(m,n) = (3.5, 4.5)$, $\bar{S}$ given by (4.2) with $\xi = 0.1, \cdots, 0.314$ and $\int_0^L h_{\text{eq}} s_{\text{eq}}\, dx = 2$. (Left) $h_{\text{eq}}$ profiles; (middle) $s_{\text{eq}}$ profiles; (right) plots of min $h_{\text{eq}}$ and max $s_{\text{eq}}$ vs the shift parameter $\xi$.

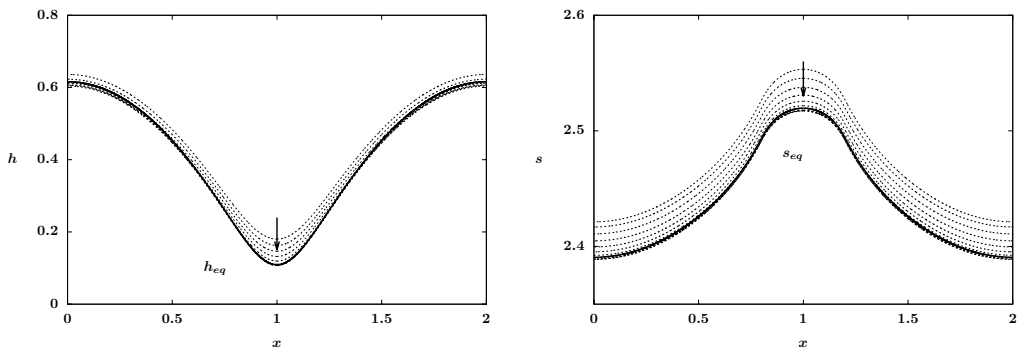

Figure 6: Convergence of PDE solutions $h$ and $s$ (in dashed lines) to the equilibrium plotted in solid lines which satisfy the ODEs (4.3) with $(m,n) = (3.5, 4.5)$ and $\bar{S}$ from (4.2) with $\xi = 0.2$.

The dynamics in Fig. 7 can be understood by looking at Fig. 8 where the minimum film thickness of the equilibrium $h_{eq}$ is plotted in terms of the shift coefficient $\xi$ in (4.2) and mobility coefficient $m$. It is shown in Fig. 8 that for $0 < \xi < 0.1$, the equations (1.1) has an equilibrium $(h_{eq}, s_{eq})$ for all positive $m$, while for $\xi \geq 0.1$ there exists a critical $m_c$ such that for $m \geq m_c$ the equilibrium solution $(h_{eq}, s_{eq})$ to the system (4.3) exist. Specifically, for $\xi = 0.3$ in (4.2), results in Fig. 8 indicate that equilibrium solutions that satisfy the ODE system (4.3) only exist for $m > m_c \approx 3.26$. Therefore the threshold $m = m_c$ divides the long-time behaviors of the PDE solutions into the two cases: infinite-time thinning with $m < m_c$ and convergence to equilibrium with $m > m_c$.

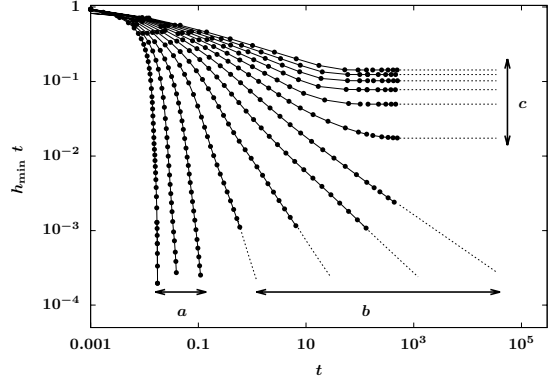

Figure 7: Plots of minimum film thickness $h_{\min}$ with identical system parameters used in Fig. 4 and $\bar{S}$ from (4.2) with $\xi = 0.3$. Region (a): $m = 0, 0.5, 1$; Region (b): $m = 1.5, 2, 2.5, 3$; Region (c): $m = 3.5, 4, \cdots, 6$. Numerical results are represented by dots, showing finite time singularity for region (a), infinite time thinning for region (b) following predictions from (4.4) plotted in dashed lines, and convergence of PDE solutions to steady states for region (c).

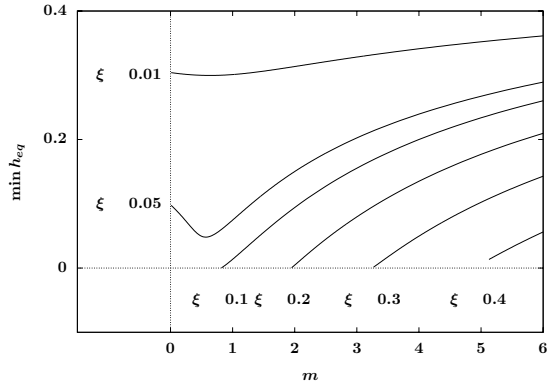

Figure 8: Plots of the minimum film thickness of the steady states $h_{eq}$ against system parameters $m$ and $n = m + 1$ with a sequence of $\xi$ values, showing that for $\xi > 0.1$ steady state solutions to (4.3) cease to exist when $m$ is smaller than a critical value $m_c$. In particular, for $\xi = 0.3$ the critical value of $m$ is $m_c = 3.26$.

## 5. Conclusions

In this work, the proof of global and local existence of strong solutions to the generalized tear film rupture model (1.1) with different system parameters regimes has been carried out. More precisely, we have shown that with mobility exponents $n = m+1$ and $3 \leq m \leq 4$ strong solutions to (1.1) exist globally, and local strong solutions to the model exist for the regime $m \geq 0$, $n \geq 0$.

The numerical results in section 4 support the conclusion of the regularity and existence of solutions in section 3. Specifically, the long time behavior of the PDE solutions to (1.1) with a family of $\bar{S}$ profiles (4.2) is investigated. For the case $n = m + 1$ and $3 \leq m \leq 4$, if an equilibrium solution $(h_{\text{eq}}, s_{\text{eq}})$ can be established in the PDE system associated with a specified total mass of salt $Q$, the PDE solutions approach to the equilibrium solution in the long time. Otherwise, without the attraction of the equilibrium, infinite-time non-uniform thinning in $h$ is expected to happen. While with mobility exponents $(m, n)$ outside the above region, for instance, with $(m, n) = (0, 3)$ in the physical model (2.10), we numerically capture the formation of finite time singularity in both $h$ and $s$ driven by the non-conservative terms in the model.

Several interesting questions regarding the PDE (1.1) remain to be solved. First, in this paper we have restricted our attention to the scenario where $n = m+1$ and $3 \leq m \leq 4$ for the proof of global existence of solutions. However, as is suggested by the sequence of simulations shown in Fig. 4 and 7, the parameter range for the existence of strong solutions can possibly be extended to larger regions. Inspired by the convergence of PDE solutions to equilibrium solutions in some of the numerical simulations, we are also interested to study whether the equilibrium solutions to (1.1) are all global attractors. More specifically, we may ask: does any solution converge to the equilibrium solutions, or how the existence of those equilibrium solutions depend on $\bar{S}$, and system parameters $(m, n)$.

**Acknowledgment**

This work was supported by the National Science Foundation under Grant No. DMS-1514826 and KI-Net RNMS11-07444.